\documentclass{amsart}

\newtheorem{theorem}{Theorem}[section]
\newtheorem{lemma}[theorem]{Lemma}

\theoremstyle{definition}
\newtheorem{definition}[theorem]{Definition}

\theoremstyle{remark}
\newtheorem{remark}[theorem]{Remark}

\newtheorem{claim}{Claim}

\begin{document}

\author{Tomoo Yokoyama}

\email{yokoyama@math.sci.hokudai.ac.jp}

\title{Codimension one minimal foliations whose leaves have fundamental groups with the same polynomial growth}

\address{
Creative Research Institution,  
Hokkaido University,  
North 21, West 10, Kita-ku, Sapporo 001-0021, Japan}

\thanks{
The author is partially supported 
by the JST CREST Program at Creative Research Institution,  
Hokkaido University.}

\maketitle

\begin{abstract}
Let $\mathcal{F}$ be a transversely orientable 
codimension one minimal foliation without vanishing cycles 
of a 
manifold $M$ and $k \in \mathbb{Z}_{\geq 0}$.  
We show that 
if the fundamental group of each leaf of $\mathcal{F}$ 
has polynomial growth of degree $k$,  
then the foliation $\mathcal{F}$ is without holonomy.   
\end{abstract}

\section{Introduction} 
In \cite{YT},  we have considered the following question: 
how is the minimal foliation without null homotopic closed transversals if the fundamental group of
each leaf is isomorphic to an elementary group? 

We have showed the following result in the paper: 
If the fundamental groups of all leaves are isomorphic to $\mathbb{Z}$, then 
the foliation is without holonomy. 
But if the fundamental groups of the leaves are ``complicated'', then  
the foliation can have nontrivial holonomy. 
In particular, there is such a foliation without vanishing cycles (e.g. Example 3.2. \cite{YT}). 
In this paper, we consider the above question for the minimal foliations without vanishing cycles 
and obtain the following result.

\begin{theorem}\label{th:001}
Let $\mathcal{F}$ a codimension one transversely orientable $C^0$
minimal foliation without vanishing cycles of a manifold $M$ and 
 $k$ a nonnegative integer. 
If 
all fundamental groups of every leaves of $\mathcal{F}$ have the same polynomial growth of degree $k$,  
then the foliation $\mathcal{F}$ is without holonomy.  
\end{theorem}
 
Note that 
a finitely generated group has polynomial growth if and only if it is virtually nilpotent \cite{G}.  
Notice that we do not require that the manifolds are compact.

\section{Preliminaries}

First, we recall a vanishing cycle in the sense of Novikov \cite{N}. 

\begin{definition} 
Let $\mathcal{F}$ be a foliation of a manifold $M$. 
A loop  $\gamma$ on the leaf $L \in F$ is called a vanishing cycle 
if there is a mapping $F: S^1 \times [0,1] \to M$ such that 
arcs $F( x, [0,1])$ for every $x \in S^1$ are transverse to $\mathcal{F}$, 
each loop $F( S^1, t)$ for any $t \in [0,1]$ is contained in some leaf 
$L_t$,   
$[F( S^1 ,0)] = [\gamma] \neq 1 \in \pi_1(L)$, and 
$[F( S^1 ,s)]  = 1 \in \pi_1(L_s)$ for all $s \in ]0,1]$. 
\end{definition}

We define a trivial fence as follows  \cite{Y}.

\begin{definition}[Trivial fences] 
Let $\mathcal{F}$ be a transversely oriented foliation of a manifold $M$. 
For a compact subset $K$ of a leaf $L_0$ of $\mathcal{F}$, 
an embedding $F: K \times [0, 1 ] \longrightarrow M$ is called a positive trivial  fence if 
$F( K \times \{ t \} ) $ is contained in a leaf $L_t$ of $\mathcal{F}$, 
$F|{K \times \{ 0 \} } $ is the inclusion $K \subset L_0$ and 
$F( \{ x \} \times [0, 1 ] ) $ is transverse to $\mathcal{F}$. 
\end{definition}

\begin{remark}
For an arcwise connected compact subset  $K$ of a leaf $L$, 
a positive trivial  fence on $K$ exists if and only if 
the holonomy on the positive side along any loop in $K$ is trivial.  
\end{remark}

We define polynomial growth for finitely generated groups.

\begin{definition}
Let $G$ a finitely generated group and 
$\mathcal{G}$ be a generating set of $G$. 
The word length function $f_G: G \to \mathbb{Z}_{\geq 0}$ on $G$ 
with respect to $\mathcal{G}$ is given by $f_G(e) := 0$ and 
$f_G( g) := \min \{ n \in \mathbb{Z}_{> 0} \mid g = g_1 \cdots g_n 
\textrm{ for some } g_i \in \mathcal{G} \cup \mathcal{G}^{-1} \}$ 
for any nontrivial element $g \in G$. 
Let $G_m$ be the set of the elements of $G$ whose word lengths are at
most $m$. 
We say that 
$G$ has polynomial growth of degree at least some nonnegative integer $k$ if 
there is a positive constant $\alpha \in  \mathbb{R}_{> 0}$ such that 
$\# G_m \geq \alpha m^k$ for all $m \in  \mathbb{Z}_{\geq 0}$. 
If the group $G$ does not have polynomial growth of degree at least $k + 1$, 
then $G$ is said to have polynomial growth of degree $k$. 
\end{definition}

Note that the definition of polynomial growth is independent of the choice of generating sets. 
To prove the main theorem, the following statement will be helpful. 

\begin{lemma}
\label{lem:006}
Let $G$ be a group with polynomial growth of degree $k \in \mathbb{Z}_{\geq 0}$ and 
$H \leq G$ subgroup of $G$ with the same growth.  
Then, for any $g \in G$, there is a positive integer $n \in \mathbb{Z}_{> 0}$ such that 
$g^n \in H$. 
\end{lemma}

\begin{proof}
Let $\mathcal{G}$ be a generating set of $G$ and   
$\mathcal{H}$ a generating set of $H$ with $\mathcal{H} \subseteq \mathcal{G}$. 
We define the word lengths on $G$ and $H$ by these generating sets. 
Let $G_m$ be the set of the elements of $G$ whose word lengths are at
 most $m$ and 
$H_m$ the set of the elements of $H$ whose word lengths are at most $m$. 
By contradiction, 
suppose that there exists $g_0 \in G$ such that 
$g_0^n \notin H$ for any $n \in \mathbb{Z}_{> 0}$. 
With out loss of generality, 
we may assume that 
$g_0 \in \mathcal{G}$. 

\begin{claim}
For $n,m \in \mathbb{Z}$ and $h, h' \in H$,  
$g^n_0 h' = g^m_0 h$ if and only if 
$n=m$ and $h'= h$. 
In particular,  $g^j_0 H_k \cap g^{j'}_0H_{k'} = \emptyset$ for any $j \neq j'$. 
\end{claim}

Indeed, 
if $n=m$ and $h'= h$, then  obviously $g^n_0 h' = g^m_0 h$. 
Conversely, 
suppose that $g^n_0 h' = g^m_0 h$.   
If $n \neq m$, then $ H \ni h'h^{-1} = g_0^{m-n} \notin H$ which is a contradiction. 
Thus $n=m$ and so $h' =h$. 

Since $G_m \supseteq \bigsqcup_{j=0}^{m} g_0^j H_{m-j}$, we have 
$$  \# G_m \geq \sum_{j=0}^m  \# g_0^j H_{m-j}.  $$
By the definition of polynomial growth, there exists $\alpha \in \mathbb{R}_{>0}$ such that 
$ \# H_m \geq \alpha m^k$ for all $m \in \mathbb{Z}_{\geq 0}$. 
Thus     
$$ \sum_{j=0}^m  \# g_0^j H_{m-j} =  \sum_{j=0}^m \# H_{j} \geq  \alpha  \cdot \sum_{j=0}^m j^k.$$
Therefore $\# G_m  \geq  \alpha \cdot  \sum_{j=0}^m j^k$.
The well known fact that  
\(
\sum_{j = 1}^{M} j^k \approx \frac{1}{k+1} M^{k+1}
\)
for sufficiently large $M \gg k$, implies that 
there is a positive constant $c \in \mathbb{R}_{>0}$ such that  
$\sum_{j=0}^m j^{k} \geq  c \cdot m^{k+1}$ for all $m \in \mathbb{Z}_{\geq 0}$.
Then we have 
$\# G_m  \geq  \alpha c \cdot m^{k+1}$ for all $m \in \mathbb{Z}_{\geq 0}$.
Therefore the growth of $G$ is more than $k$. 
This contradicts with the hypothesis.
\end{proof}

\section{A Key lemma and the proof of the main theorem}

The following key lemma will complete the proof of main theorem. 
Denote by $\mathrm{gr}(G)$ the degree of the polynomial growth of a group $G$. 

\begin{lemma}
\label{prop:005}
Let $\mathcal{F}$ a codimension one transversely orientable  
minimal $C^0$ foliation on a 
manifold $M$. 
If there is a leaf $L_0 \in \mathcal{F}$ without holonomy such that 
the induced homomorphism $i_*: \pi_1(L_0) \to \pi_1(M)$ of the inclusion $i: L_0 \to M$ is injective 
and that 
$\mathrm{gr}(\pi_1 (L_0)) \geq \mathrm{gr}(\pi_1(L))$ for any leaf $L \in \mathcal{F}$,  
then $\mathcal{F}$ is without holonomy. 
\end{lemma}

\begin{proof}
Take a compact subset $K := \bigvee_k S^1 \subset L_0$ with a base point $x$ for the root of this bouque 
such that 
the homomorphism $\pi_1 (K, x) \to \pi_1 (L_0, x)$ induced by the inclusion $K \to L_0$ is surjective.  
Let $F: K \times [0,1] \to M$ be a positive trivial fence on the leaf $L_0$, 
$i_{t}: L_t \to M$ inclusions, 
$F_t : (K, x) \to (L_t, x_t)$ the induced maps by $i_{t} \circ F_t = F(\cdot, t)$, 
and 
$H_t := \mathrm{im}(F_{t*}) \leq \pi_1(L_t, x_t)$ subgroups for any $t \in [0,1]$.  
By the definitions and the injectivity,  
$\mathrm{gr}(\pi_1(L_0, x_0)) 
= \mathrm{gr}(H_0)
= \mathrm{gr}(\mathrm{im}(F_{0*})) 
= \mathrm{gr}(i_{0*}(\mathrm{im}(F_{0*}))) 
= \mathrm{gr}(\mathrm{im}((i_0 \circ F_{0})_{*}))  $.
Since 
$\mathrm{im}(i_0 \circ F_{0}) = F(K, 0) $ and $\mathrm{im}(i_t \circ F_{t}) = F(K, t)$ 
are free homotopic in $M$, 
we have 
$\mathrm{gr}(\mathrm{im}((i_0 \circ F_{0})_{*})) 
= \mathrm{gr}(\mathrm{im}((i_t \circ F_{t})_{*})) 
= \mathrm{gr}(i_{t*}(\mathrm{im}(F_{t*}))) 
= \mathrm{gr}(i_{t*}(H_t)) 
\leq \mathrm{gr}(H_t) \leq \mathrm{gr}(\pi_1 (L_t, x_t))$ for any $t \in [0,1]$. 
By the hypothese,  we obtain 
$\mathrm{gr}(\pi_1 (L_t, x_t)) \leq \mathrm{gr}(\pi_1(L_0, x_0))$ for any $t \in [0,1]$. 
Thus 
$ H_t$ and $\pi_1 (L_t)$ have the same polynomial growth.   
By Lemma \ref{lem:006}, 
we have that 
for any $t \in [0, 1]$ and any $g \in \pi_1(L_t)$, 
there exists some $n \in \mathbb{Z}_{> 0}$ such that  
$g^n \in \mathrm{im}F_{t*}$. 
Since any codimension one transversely orientable foliation has no finite holonomy, 
the leaves in the saturation of $F( K \times ]0, 1[)$ are without holonomy.  
By the minimality of $\mathcal{F}$, 
we obtain 
the saturation is $M$. 
Thus $\mathcal{F}$ is without holonomy.  
\end{proof}

Now we prove the main theorem. 

\begin{proof}[
Proof of Theorem\ref{th:001}] 
Since 
the union of the leaves without holonomy is dense G$_{\delta}$ \cite{EMT}, 
there is a leaf $ L$ of $\mathcal{F}$ without holonomy. 
Since  $\mathcal{F}$ has no vanishing cycles, 
by Theorem 3.4. p.147 \cite{HH},  
we obtain that 
the homomorphism $\pi_1(L) \to \pi_1(M)$ induced by the inclusion 
$L \to M$  is injective. 
By Lemma \ref{prop:005},  
we have 
$\mathcal{F}$ is without holonomy. 
\end{proof}

\section{Remarks}

The condition that all fundamental groups of leaves have the same polynomial growth 
is necessary 
because Anosov foliations for flows satisfy the conditions except this growth condition  
and 
have nontrivial holonomy. 

The following example shows that 
the minimal condition in the main theorem is necessary. 
Consider a continuous map $f : [0, 1] \to [0, 1]$ such that $f( 0) = 0, f(1) = 1, $ and $ f(x ) > x $ for any $x \in ]0, 1[$. 
Identifying $0$ with 1, we obtain the map induced on $S^1$ by $f$. 
Take a suspension foliation $\mathcal{F}_1$ on $T^2$ of the induced map. 
Then $\mathcal{F}_2 = \{ L \times \mathbb{S}^1 \mid L \in \mathcal{F}_1 \}$ is a foliation on $T^3$
which consists of one torus and other cylinders.  
Remove the simple closed curve in the toral leaf parallel to the suspended direction. 
Then the resulting cylindrical leaf has nontrivial holonomy and 
the resulting non-minimal foliation consists of only cylinders. 

However we do not know whether the minimal condition is necessary for compact manifolds.



\begin{thebibliography}{Gh9}

 \bibitem[EMT77]{EMT}
D. B. A.Epstein,   K. C. Millett, D. Tischler,  \textit{Leaves without holonomy.} 
 J. London Math. Soc. (2) 16 (1977),  no. 3,  548--552.
 \bibitem[G81]{G}
M. Gromov, 
\textit{Groups of polynomial growth and expanding maps} 
Inst. Hautes Etudes Sci. Publ. Math. 53 (1981), 53–73.
 \bibitem[HH83]{HH}
   G.Hector, U. Hirsch,  
 \textit{Introduction to the geometry of foliations. Part B.}
 Foliations of codimension one. Aspects of Mathematics, E3. Friedr. Vieweg \& Sohn, Braunschweig, 1983.
\bibitem[N65]{N}
S. P. Novikov, 
\textit{Topology of foliations} 
Trudy Moskov. Mat. Obshch. 14 (1965) 248--278.
\bibitem[Y09]{Y}
T. Yokoyama,
\textit{Codimension one minimal foliations and the higher homotopy groups of leaves} 
C. R. Acad. Sci. Paris, Ser. I 347 (2009) 655--658.
\bibitem[YT08]{YT}
T. Yokoyama, T. Tsuboi,
\textit{Codimension one foliations and the fundamental groups of leaves} 
Annales de L'Institut Fourier. (2008), vol. 58, no.2, 723--731.

\end{thebibliography}
\end{document}